\documentclass[12pt]{article}
\usepackage{latexsym}
\usepackage{amsfonts}
\begin{document}
\newcommand{\tr}{\text {tr}}
\newcommand{\Hom}{\text {Hom}}
\newcommand{\Irr}{\text {Irr}}
\newcommand{\End}{\text {End}}
\newcommand{\Comod}{\text {Comod}}
\newcommand{\Rep}{\text {Rep}}
\newcommand{\n}{\mathfrak{n}}
\newcommand{\g}{\mathfrak{g}}
\newcommand{\text}[1]{\mbox{{\rm #1}}}
\newcommand{\gd}{\delta}
\newcommand{\itms}[1]{\item[[#1]]}
\newcommand{\nin}{\in\!\!\!\!\!/}
\newcommand{\sub}{\subset}
\newcommand{\cntd}{\subseteq}
\newcommand{\go}{\omega}
\newcommand{\Pa}{P_{a^\nu,1}(U)}
\newcommand{\fx}{f(x)}
\newcommand{\fy}{f(y)}
\newcommand{\gD}{\Delta}
\newcommand{\gl}{\lambda}
\newcommand{\gL}{\Lambda}
\newcommand{\half}{\frac{1}{2}}
\newcommand{\sto}[1]{#1^{(1)}}
\newcommand{\stt}[1]{#1^{(2)}}
\newcommand{\Z}{\hbox{\sf Z\kern-0.720em\hbox{ Z}}}
\newcommand{\singcolb}[2]{\left(\begin{array}{c}#1\\#2
\end{array}\right)}
\newcommand{\ga}{\alpha}
\newcommand{\gb}{\beta}
\newcommand{\gga}{\gamma}
\newcommand{\ul}{\underline}
\newcommand{\ol}{\overline}
\newcommand{\qed}{\kern 5pt\vrule height8pt width6.5pt depth2pt}
\newcommand{\Lrraro}{\Longrightarrow}
\newcommand{\Nb}{|\!\!/}
\newcommand{\NN}{{\rm I\!N}}
\newcommand{\bsl}{\backslash}
\newcommand{\gt}{\theta}
\newcommand{\op}{\oplus}
\newcommand{\C}{k}
\newcommand{\Q}{{\bf Q}}
\newcommand{\Op}{\bigoplus}
\newcommand{\CR}{{\cal R}}
\newcommand{\grr}{\omega_1}
\newcommand{\ben}{\begin{enumerate}}
\newcommand{\een}{\end{enumerate}}
\newcommand{\ndiv}{\not\mid}
\newcommand{\bab}{\bowtie}
\newcommand{\hal}{\leftharpoonup}
\newcommand{\har}{\rightharpoonup}
\newcommand{\ot}{\otimes}
\newcommand{\OT}{\bigotimes}
\newcommand{\bwe}{\bigwedge}
\newcommand{\gep}{\varepsilon}
\newcommand{\rbraces}[1]{\left( #1 \right)}
\newcommand{\bbox}{$\;\;\rule{2mm}{2mm}$}
\newcommand{\sbraces}[1]{\left[ #1 \right]}
\newcommand{\bbraces}[1]{\left\{ #1 \right\}}
\newcommand{\OO}{_{(1)}}
\newcommand{\TT}{_{(2)}}
\newcommand{\FF}{_{(3)}}
\newcommand{\minus}{^{-1}}
\newcommand{\CV}{\cal V}
\newcommand{\CVs}{\cal{V}_s}
\newcommand{\un}{U_q(sl_n)'}
\newcommand{\on}{O_q(SL_n)'}
\newcommand{\slq}{U_q(sl_2)}
\newcommand{\olq}{O_q(SL_2)}
\newcommand{\UU}{U_{(N,\nu,\go)}}
\newcommand{\HH}{H_{n,q,N,\nu}}
\newcommand{\ct}{\centerline}
\newcommand{\bs}{\bigskip}
\newcommand{\qua}{\rm quasitriangular}
\newcommand{\ms}{\medskip}
\newcommand{\noin}{\noindent}
\newcommand{\mat}[1]{$\;{#1}\;$}
\newcommand{\raro}{\rightarrow}
\newcommand{\map}[3]{{#1}\::\:{#2}\raro{#3}}
\newcommand{\alg}{{\rm Alg}}
\def\newtheorems{\newtheorem{theorem}{Theorem}[subsection]
                 \newtheorem{cor}[theorem]{Corollary}
                 \newtheorem{prop}[theorem]{Proposition}
                 \newtheorem{lemma}[theorem]{Lemma}
                 \newtheorem{defn}[theorem]{Definition}
                 \newtheorem{Theorem}{Theorem}[section]
                 \newtheorem{Corollary}[Theorem]{Corollary}
                 \newtheorem{Proposition}[Theorem]{Proposition}
                 \newtheorem{Lemma}[Theorem]{Lemma}
                 \newtheorem{Defn}[Theorem]{Definition}
                 \newtheorem{Example}[Theorem]{Example}
                 \newtheorem{Remark}[Theorem]{Remark}
                 \newtheorem{claim}[theorem]{Claim}
                 \newtheorem{sublemma}[theorem]{Sublemma}
                 \newtheorem{example}[theorem]{Example}
                 \newtheorem{remark}[theorem]{Remark}
                 \newtheorem{question}[theorem]{Question} 
                 \newtheorem{Definition}[Theorem]{Definition}              
                 \newtheorem{Question}[Theorem]{Question}
                 \newtheorem{conjecture}{Conjecture}[subsection]
                 \newtheorem{Conjecture}[Theorem]{Conjecture}
}
\newtheorems
\newcommand{\proof}{\par\noindent{\bf Proof:}\quad}
\newcommand{\dmatr}[2]{\left(\begin{array}{c}{#1}\\
                            {#2}\end{array}\right)}
\newcommand{\doubcolb}[4]{\left(\begin{array}{cc}#1&#2\\
#3&#4\end{array}\right)}
\newcommand{\qmatrl}[4]{\left(\begin{array}{ll}{#1}&{#2}\\
                            {#3}&{#4}\end{array}\right)}
\newcommand{\qmatrc}[4]{\left(\begin{array}{cc}{#1}&{#2}\\
                            {#3}&{#4}\end{array}\right)}
\newcommand{\qmatrr}[4]{\left(\begin{array}{rr}{#1}&{#2}\\
                            {#3}&{#4}\end{array}\right)}
\newcommand{\smatr}[2]{\left(\begin{array}{c}{#1}\\
                            \vdots\\{#2}\end{array}\right)}

\newcommand{\ddet}[2]{\left[\begin{array}{c}{#1}\\
                           {#2}\end{array}\right]}
\newcommand{\qdetl}[4]{\left[\begin{array}{ll}{#1}&{#2}\\
                           {#3}&{#4}\end{array}\right]}
\newcommand{\qdetc}[4]{\left[\begin{array}{cc}{#1}&{#2}\\
                           {#3}&{#4}\end{array}\right]}
\newcommand{\qdetr}[4]{\left[\begin{array}{rr}{#1}&{#2}\\
                           {#3}&{#4}\end{array}\right]}

\newcommand{\qbracl}[4]{\left\{\begin{array}{ll}{#1}&{#2}\\
                           {#3}&{#4}\end{array}\right.}
\newcommand{\qbracr}[4]{\left.\begin{array}{ll}{#1}&{#2}\\
                           {#3}&{#4}\end{array}\right\}}

\title{On Cotriangular Hopf Algebras}
\author{Pavel Etingof
\\Department of Mathematics, Rm 2-165\\
Massachusetts Institute of Technology\\
77 Massachusetts Avenue\\
Cambridge, MA 02139
\and Shlomo Gelaki\\Mathematical Sciences Research Institute\\
1000 Centennial Drive\\Berkeley, CA 94720}
\maketitle

\section{Introduction}
In [EG1, Theorem 2.1] we proved that any triangular
semisimple Hopf
algebra over an algebraically closed field $\C$ of characteristic $0$
is obtained from the group algebra $\C[G]$ of a finite group $G,$
by twisting its  comultiplication by a twist in the sense of
Drinfeld [Dr]. Since
semisimple Hopf algebras are finite-dimensional, dualizing yields
that any cotriangular semisimple Hopf algebra over $\C$ is
obtained from $\C[G]^*,$ the function algebra on $G,$ by
twisting its  multiplication by a Hopf $2-$cocycle in the
sense of Doi
[Do] (see Section 2 below).

In this paper we generalize Theorem 2.1 from [EG1] to
not necessarily finite-dimensional cotriangular Hopf
algebras $A$ over $\C.$ Namely, our main result (see Section 3
below) is:

\noin
{\bf Theorem} A cotriangular Hopf algebra $A$ over $k$ is obtained
from the
function algebra ${\cal O}(G)$ of a pro-algebraic
group $G,$ by twisting its  multiplication by a Hopf
$2-$cocycle, and possibly changing its R-form by a central grouplike
element of $A^*$ of order $\le 2,$ if and only if
$\tr (S^2|_C)=\dim(C)$ 
for any finite-dimensional subcoalgebra $C$ of $A$ (where $S$ is
the antipode
of $A$). 

The main challenge in the proof of this theorem 
(see Section 4 below) is to establish 
the ``if'' direction. 
The key step in the proof of the ``if'' part  is
to show that 
our trace condition on $A$ guarantees that 
the categorical dimensions of
objects in the category of its finite-dimensional right comodules are
non-negative integers (maybe after modifying the R-form). 
This enables us to apply the same theorem of 
Deligne on Tannakian categories [De]
that we applied in the proof of Theorem 2.1 from [EG1].

In Section 5, we give examples of twisted function algebras. 
In particular, we show that 
in the infinite-dimensional case, the 
squared antipode for such an algebra may not equal the identity 
(see Example \ref {2} below).

In Section 6, we show that in all of our examples, 
the operator $S^2$ is unipotent on $A$, and conjecture it to 
be the case for any twisted function algebra. 
We prove this conjecture, using the quantization theory of [EK1-2],
in a large number of special cases. 

In Section 7, we formulate a few open  
questions.

Throughout the paper, $k$ will denote an algebraically closed field of
characteristic $0.$

\noin
{\bf Acknowledgements} The first author is grateful to Ben Gurion
University for its warm hospitality, and to Miriam Cohen and the
Dozor Fund for making his visit possible; 
his work was also supported by the NSF grant DMS-9700477.

The second author is grateful to Susan
Montgomery for numerous useful conversations.

The authors would like to acknowledge that this paper was inspired
by the work [BFM].   

\section{Hopf $2-$cocycles}

Let $A$ be a coassociative coalgebra over $k.$
For $a\in A,$ we write $\Delta(a)=\sum a_1\ot a_2,$ \linebreak
$(I\ot \Delta)\Delta(a)=\sum a_1\ot a_2\ot a_3$ etc, where $I$ denotes
the identity map of $A.$

Recall that $A^*$ is an associative algebra with product defined by 
$(f*g)(a)=\sum f(a_1)g(a_2).$
This product is called the {\em convolution product}. 

Now let $(A,m,1,\Delta, \varepsilon,S)$ be a Hopf algebra over
$k.$ 

Recall [Do] that a linear form $J:A\ot A\raro k$ is called a {\em
Hopf
$2-$cocycle} for $A$ if it has an inverse $J^{-1}$ under the convolution
product $*$ in
$\Hom _k (A\ot A,k)$, and satisfies:
\begin{equation}\label{2coc}
\sum J (a_1b_1,c)J (a_2,b_2)= \sum J
(a,b_1c_1)J (b_2,c_2)\;\text {and}\; J (a,1)=\varepsilon(a)=J (1,a)
\end{equation}
for all $a,b,c\in A.$ 

Given a Hopf $2-$cocycle $J$ for $A,$ one can construct a new Hopf
algebra $(A^{J}, m^{J},1, \Delta, \varepsilon,S^{J})$ as
follows. As a coalgebra, $A^{J}=A.$ The new multiplication is given by 
\begin{equation}\label{nm}
m^{J}(a\ot b)=\sum J^{-1} (a_1,b_1)a_2b_2J 
(a_3,b_3)
\end{equation}
for all $a,b\in A.$ The new antipode is given by
\begin{equation}\label{na}
S^{J}(a)=\sum J ^{-1}(a_1,S(a_2))S(a_3)J 
(S(a_4),a_5)
\end{equation}
for all $a\in A.$

Suppose $A$ is also co(quasi)triangular with
universal R-form $R:A\ot A\raro k$ (see e.g. [K,
VIII.5.1]). Then it is straightforward to verify
that $A^{J}$ is co(quasi)triangular with universal R-form
$R^{J}:A^{J}\ot A^{J}\raro k,$ where
$R^{J}:=(J\circ \tau)^{-1}*R*J$ (here
$\tau:A\ot A\raro A\ot A$ is the usual flip map).

Recall [Dr] that a {\em twist} for a Hopf algebra $B$ is an
invertible element $J\in B\ot B$ which satisfies
\begin{equation}\label{t1}
(\Delta\ot I)(J)(J\ot 1)=(I\ot \Delta)(J)(1\ot J)\;
\text{and}\;(\varepsilon\ot
I)(J)=(I\ot \varepsilon)(J)=1.
\end{equation}
It is straightforward to verify that if $A$ is finite-dimensional, then
$J\in A^*\ot A^*$ is a Hopf $2-$cocycle for $A$ if and only if it is a
twist for $A^*.$

\section{The Main Theorem} 
Let $(A,R)$ be a cotriangular Hopf algebra
over $\C$ (not necessarily
finite-dimensional). Define the Drinfeld element of $(A,R)$ to be
the linear form $u:A\raro \C$ determined by
\begin{equation}\label{drinf}
u(a)=\sum R(a_2,S(a_1)).
\end{equation}
Recall that $u\in A^*$ is a grouplike element, and that
\begin{equation}\label{antip}
S^2(a)=(u*I*u^{-1})(a)=\sum u(a_1)a_2u^{-1}(a_3)
\end{equation}
for all $a\in A.$

Suppose $c\in A^*$ is a central grouplike element of order $\le 2,$ and set
\begin{equation}\label{rc}
R_c:=\frac{1}{2}(\varepsilon \ot \varepsilon +\varepsilon \ot c+c\ot
\varepsilon -c\ot c).
\end{equation}
Then it is straightforward to verify that $(A,R*R_c)$ is
cotriangular with
Drinfeld element $u*c.$

Note that by (\ref{antip}), $S^2$ preserves any subcoalgebra of $A.$
\begin{Definition}\label{pseu}
We say that $A$ is pseudoinvolutive if
$\tr(S^2|_{C})=\dim(C)$ for
any finite-dimensional subcoalgebra $C$ of $A.$
\end{Definition}
\begin{Remark}\label{fd}
{\rm If $A$ is finite-dimensional, then
pseudoinvolutivity is equivalent to involutivity ($S^2=I$).
Indeed, by [R], 
$S^2$ has a finite order, so its eigenvalues on $A$ are roots of
$1.$ But a sum of
$\dim(A)$ roots of $1$ can be equal to $\dim(A)$ only if all
of them are
$1.$} $\Box$
\end{Remark}

We can now state our main theorem.
\begin{Theorem}\label{equiv}
Let $(A,R)$ be a cotriangular Hopf algebra over $\C$.
 Then the following two conditions are equivalent:

\noin
(i) $A$ is pseudoinvolutive.

\noin
(ii) $A$ is a twisted function algebra ${\cal 
O}(G)^{J}$ on a
pro-algebraic group $G,$ and furthermore, there exists a
central grouplike element $c\in {\cal O}(G)^*$ of order $\le
2,$ such that $(A,R)$ is isomorphic to $({\cal
O}(G)^{J},(J\circ
\tau)^{-1}*R_c*J)$ as
cotriangular Hopf algebras.
\end{Theorem}
\begin{Remark}\label{gen}
{\rm Theorem \ref{equiv} is a generalization of Theorem 2.1 from
[EG1]. Indeed, let $A$ be a finite-dimensional triangular Hopf algebra
over $k.$ Equivalently, $A^*$ is a finite-dimensional cotriangular
Hopf
algebra over $k.$ Now, by Remark \ref{fd}, $A^*$ is
pseudoinvolutive if and
only if $S^2=I.$ By [LR], this is equivalent to the semisimplicity
of $A$ and
$A^*.$ Hence by Theorem \ref{equiv}, $A^*$ is a finite-dimensional
semisimple
cotriangular Hopf algebra over $k$ if and only if it is a
twisted function
algebra ${\cal O}(G)^{J}$ on a pro-algebraic group
$G$ (possibly, with a changed cotriangular structure). But of course,
$G$ must be a finite group.} $\Box$
\end{Remark}
\begin{Remark}\label{unique}
{\rm The data $(G,c,J)$ corresponding to a pseudoinvolutive
cotriangular Hopf algebra over $k,$ is unique up to isomorphism of
such triples and gauge transformations of $J$
(see [EG2]). The proof is similar to the proof of Lemma 3.5
in [EG2].} $\Box$
\end{Remark}
\begin{Remark}\label{susan}
{\rm If $(A,R)$ is a cocommutative cotriangular Hopf
algebra over $k$ (hence also commutative), then Theorem
\ref{equiv} is applicable since in this situation $S^2=I.$
Thus, [BFM, Theorem 3.19(i)],
which claims that in this case $(A,R)$ is a twisted 
group algebra (maybe with $R\to R*R_c$), is a special 
case of our result.} $\Box$
\end{Remark}
\begin{Remark}\label{miriam}
{\rm In [CWZ, Theorem 2.1] the authors prove that if $(A,R)$ is
cotriangular and its Drinfeld element $u$ acts as the identity on a
finite-dimensional right $A$-comodule $V,$ then the
characters of the usual action of the symmetric group $S_n$ on
$V^{\ot n}$ and
the one arising from the braiding $R$ are equal. We note that
this result is a consequence of Theorem \ref{equiv}. Namely, one
should apply the theorem to the cotriangular Hopf algebra $H_V\subset
A$ which is generated by the elements of the form $(f\otimes
I)(\rho_{_V}(v)),$ $v\in V,f\in V^*$ (here $\rho_{_V}$ denotes the
structure map of $V$), and get that the symmetric category of
finite-dimensional right comodules of the cotriangular Hopf algebra
$H_V$ is equivalent
to that of ${\cal O}(G)$ for some pro-algebraic group $G$.
Since the characters of $S_n$ on $V^{\ot n}$ are invariant under
equivalences of symmetric
categories, the result follows. We also see that Theorem 2.1 of [CWZ]
can be strengthened by replacing the assumption $u|_V=1$ by the weaker
assumption that $u|_V$ is unipotent.} $\Box$
\end{Remark}
\section{The Proof of Theorem \ref{equiv}} 
\subsection{Changing the Cotriangular Structure}
In this subsection we prove the following proposition, which is one
of the key ingredients in the proof of Theorem \ref{equiv}.  
\begin{prop}\label{changing}
For any pseudoinvolutive
cotriangular Hopf algebra $(A,R)$ over $k,$ there exists a 
central grouplike element $c\in A^*$ of order $\le 2$ such that 
for the cotriangular Hopf algebra $(A,R*R_c)$, with Drinfeld element
$u,$ one has $\text{tr}(u|_V)=\text{dim}(V)$
for any finite-dimensional right $A$-comodule $V.$
\end{prop}

The rest of the subsection is devoted to the proof of this
proposition. 

Let $(A,R)$ be a cotriangular Hopf algebra
over $k.$ Let ${\cal C}:=\Comod _{f.d}(A)$ be the category of
finite-dimensional right
$A-$comodules, and $\Irr({\cal C})$ be the 
subcategory of all irreducible
objects of ${\cal C}.$ It is straightforward to check that ${\cal C}$ is
an abelian rigid symmetric tensor category in the sense of [DM]. The
unit object of ${\cal C}$ is ${\bf 1}:=k,$ and clearly
$\End({\bf 1})=k.$ Also, recall that for any object $V\in 
{\cal C},$ one can define its {\em categorical dimension} [DM],
denoted by $\dim_c(V),$ to be the image of $1$ under the morphism
$\C\raro V\ot V^*\raro V^*\ot V\raro \C$ (where the morphism $V\ot
V^*\raro V^*\ot V$ is the braiding map). It is straightforward to
verify that
\begin{equation}\label{dimc}
\dim_c(V)=\tr|_V(u),
\end{equation}
where $u$ is regarded as the linear
map $V\raro V$ determined by $v\mapsto (I\ot u)\rho_{_V}(v)$
(where $\rho_{_V}$ denotes the structure map of $V$). Observe that
$\tr|_V(u)=\tr|_V(u^{-1})$ for any $V\in {\cal C}.$ Indeed, we have
that $u^{-1}=u\circ S$ and $R=R\circ (S\ot S),$ hence $u(a)=u^{-1}(a)$
for any cocommutative element $a\in A.$ But, $\tr|_V\in A$ is a
cocommutative element.

For any object $V\in {\cal C},$ set 
\begin{equation}\label{av}
A_V:=\{(f\ot I)\rho_{_V}(v)|v\in V,\,f\in V^*\}.
\end{equation}
It is clear
that $A_V$ is a finite-dimensional subcoalgebra
of $A.$

 From now on we assume that $A$ is pseudoinvolutive. 

\begin{defn}\label{pn}
We say that an object $V\in {\cal C}$ is positive if
$\dim_c(V)=\dim(V),$ and negative if $\dim_c(V)=-\dim(V).$
\end{defn}

\begin{lemma}\label{simp}
An object of $\Irr({\cal C})$ is either positive or
negative, and any object of ${\cal C}$ is positive (resp. negative) if
and only if so are all the composition factors of its Jordan-H\"older
series. 
\end{lemma}
\proof Let $X\in \Irr({\cal C}).$ Then $A_X=X\otimes X^*.$ 
Since $\tr|_X(u)=\tr|_X(u^{-1}),$ it follows from
(\ref{antip}) that $\dim(X)^2=\tr
(S^2|_{A_X})= (\tr |_X(u))(\tr
|_X(u^{-1}))=\dim_c(X)^2.$ Thus, $\dim_c(X)=\pm \dim(X)$ as desired. 
Moreover, for any $V\in {\cal C},$ if $0=V_0\subset V_1\subset
\cdots\subset V_n=V$ is its Jordan-H\"older series, then
$\dim_c(V)=\sum_{i=1}^n \dim_c(V_{i}/V_{i-1}),$ and the
result follows. \qed

Consider now the abelian category of finite-dimensional right
{\em bicomodules} over $A,$ i.e. right comodules over $A\ot A^{cop}.$
It is clear that irreducible objects of this category are of the form
$X\ot
Y^*,$ where $X,Y\in \Irr({\cal C}).$

For any $V\in {\cal C},$ it is clear that
$A_V$ is a right $A-$bicomodule (recall that $A_V^*$ is the image of
$A^*\raro \End(V),$ so it is a left $A^*\ot A^{*op}$-module).
For any $X,Y\in \Irr({\cal C}),$ let $N_V(X,Y)$ be the
multiplicity of occurrence of $X\otimes Y^*$ as a composition factor
in the Jordan-H\"older series of $A_V$ regarded as a right bicomodule
over $A.$ 
\begin{lemma}\label{forbid}
For any $V\in {\cal C}$ and $X,Y\in \Irr({\cal C})$ with opposite
signs, $N_V(X,Y)=0.$
\end{lemma}
\proof Indeed, 
\begin{eqnarray*}
\lefteqn{\tr (S^2|_{A_V})}\\
&=& \sum_{X,Y\in \Irr({\cal C})} N_V(X,Y)\tr|_{X\ot
Y^*}\left (u\ot (u^{-1})^*\right)\\
&=& \sum_{X,Y\in \Irr({\cal C})} N_V(X,Y)\dim_c(X)\dim_c(Y).
\end{eqnarray*} 
But by pseudoinvolutivity, this should be equal to 
$$\dim(A_V)=\sum
_{X,Y\in \Irr({\cal C})}N_V(X,Y)\dim(X)\dim(Y).$$ 
This implies that if
$\dim_c(X)$ and $\dim_c(Y)$
have opposite signs then $N_V(X,Y)=0,$ otherwise $\sum
_{X,Y\in \Irr({\cal C})}N_V(X,Y)\dim_c(X)\dim_c(Y)< \sum
_{X,Y\in \Irr({\cal C})}N_V(X,Y)\dim(X)\dim(Y).$ \qed 
\begin{lemma}\label{ext}
For any $V,W\in {\cal C}$ with opposite signs, one has ${\rm
Ext}^1(V,W)=0.$
\end{lemma}
\proof We have to show that any exact sequence 
$0\raro V\raro U\raro W\raro 0$ in ${\cal C}$ splits.
Indeed, we have natural
coalgebra embeddings $A_V\hookrightarrow A_U,$
$A_W\hookrightarrow A_U$ (
generated by the coactions on the subcomodule and the
quotient comodule), and
the sum of them is a coalgebra embedding $A_V\oplus
A_W\hookrightarrow A_U$ (it is obvious that the sum is direct as the
coalgebras $A_V,A_W$ do not intersect, because of the opposite signs
of $V,W$). Consider the quotient
$A_U/(A_V\oplus A_W)$ as an $A-$ bicomodule. Its all composition
factors are of the form $X_+\ot X_-^*,$ where $X_+\in \Irr({\cal C})$
is positive and $X_-\in \Irr({\cal C})$ is negative, which is
impossible by Lemma \ref{forbid}. Therefore, this bicomodule must be
zero. Thus, $A_V\oplus A_W\cong A_U.$
Let $W':=\rho_U^{-1}(U\ot A_W).$ Then $U=V\oplus W',$ and we are
done. \qed
\begin{lemma}\label{eu}
Any object $V\in {\cal C}$ can be uniquely represented
as a direct sum $V_+\oplus V_-,$ where $V_+\in {\cal C}$ is
positive and $V_-\in {\cal C}$ is negative. Furthermore, for any
two objects $V,W\in {\cal C},$ $(V\ot W)_+=(V_+\ot W_+)
\oplus
(V_-\ot W_-),$ and $(V\ot
W)_-=(V_+\ot W_-) \oplus (V_-\ot W_+).$ 
\end{lemma}
\proof
We first prove the existence part of the lemma by induction in
$\dim(V).$ Let $V\in {\cal C},$ and $X\in \Irr({\cal C})$ be an
irreducible
subcomodule of $V.$ Let us assume $X$ is positive. By the induction
assumption we have $V/X=W_+\oplus W_-.$ Let $V_+$ be the preimage of
$W_+$ under the projection $V\raro W_+.$ Then $V_+$ is positive, and
$V/V_+=W_-.$ Therefore, by Lemma \ref{ext}, $V$ is isomorphic to
$V_+\oplus W_-,$ and the result follows. If $X$ is negative, the proof
is similar.

We now prove the uniqueness part of the lemma. Let $V= V_+\oplus V_-.$
Then it is easy to see that $V_+$ is the sum of all positive
subcomodules of $V,$ and $V_-$ is the sum of all negative
subcomodules, which implies uniqueness. 

The last statement of the lemma is obvious. \qed 

Let $D$ be any coalgebra over $k,$ let $\Comod(D)$ be its
category of right comodules and $F:\Comod(D)\raro Vec$ be the
forgetful functor. Recall that $\End(F)\cong D^{*}$ as
algebras. Indeed, an element $\eta\in
\End(F)$ is by definition, a collection of linear maps
$\eta_V:V\raro V,$ $V\in \Comod(D),$ which commute with
comodule morphisms. In particular, $\eta_D:D\raro D$ is a
linear map between right $D-$comodules, and it commutes with
right actions by elements of $D^*.$ Thus, $\eta_D$ comes from
a left action by an element of $D^*$
(namely, by $\eta_D^*(\varepsilon)$). Note that if moreover
$\eta_V:V\raro V$ is a comodule map for all $V\in
\Comod(D),$ then we have an endomorphism of the identity
functor $Id$ of $\Comod(D),$ and the resulting element of
$D^*$ is central. Thus, $\End(Id)\cong \text{Center}(D^*)$ as
algebras.

For any $V\in {\cal C},$ define the comodule automorphism
$c_{_V}$ of $V$ by $c_{_V}:=I$ on $V_+$ and $c_{_V}:=-I$ on
$V_-,$ where $I$ is the identity map of $V.$
\begin{lemma}\label{cv}
The collection $\{c_{_V}|V\in {\cal C}\}$ determines
a central grouplike element $c\in A^*$ of order $\le 2.$
\end{lemma}
\proof 
The collection $\{c_{_V}|V\in {\cal C}\}$ is an
element of the algebra $\End(Id),$ hence, by the
preceding remarks, determines a central element $c\in A^*.$
Now, it is clear from Lemma \ref{eu} that $c$ is a grouplike
element of order $\le 2.$ \qed

Now let us finally prove Proposition \ref{changing}.
Let $c\in A^*$ be the central grouplike element of
order $\le 2$ whose existence is guaranteed by Lemma \ref{cv}. Let
$R_c$ be
as in (\ref{rc}). Then, after changing $R$ to $R*R_{c},$ the new
Drinfeld element is $u':=u*c,$ and we get that
$\tr|_V(u')=\dim(V)$ for any object $V\in {\cal C}.$ 
The proposition is proved. \qed
\subsection{The Proof of Theorem \ref{equiv}}

$(i)\Rightarrow (ii)$.

Proposition \ref{changing} implies that without 
loss of generality, we can assume that $\text{tr}(u|_V)=
\text{dim}(V)$ for all finite-dimensional right $A-$comodules. 

Now comes the main step of the proof, which is the usage
of the following theorem of Deligne. 

\begin{theorem} [De, Theorem 7.1] Let ${\cal C}$ be an abelian
 rigid symmetric tensor category over $\C$ such that
$\End({\bf 1})=k,$ in which categorical dimensions of
objects are non-negative integers. Then there exist a
pro-algebraic group $G$ and a $k-$linear equivalence of abelian rigid
symmetric tensor categories
$F:{\cal C}\raro \Rep_{f.d}(G)$ (where $\Rep_{f.d}(G)$ is the
category of finite-dimensional algebraic $\C-$representations of $G$).
\end{theorem}

This theorem implies that in our situation, we have 
an equivalence 
$$F:\Comod_{f.d}(A)\raro \Comod_{f.d}({\cal O}(G))$$ 
of rigid symmetric tensor categories. It is obvious that $F$
preserves dimensions. 

Now we will need the following proposition, whose proof 
occupies the next subsection. 

\begin{prop}\label{morita} 
Let $A$ and $B$ be two coassociative coalgebras
with counit over $k,$ and 
$F:\Comod _{f.d}(A) \to \Comod _{f.d}(B)$ be an
equivalence 
between the abelian categories of finite-dimensional right
comodules over $A$ and $B,$ which preserves dimensions. 
Then there exists an isomorphism of coalgebras $\phi:A\to B$ 
such that $F$ is isomorphic to the direct image functor $\phi_*$.
\end{prop}
\begin{remark}
{\rm In the case when $A$,$B$ are cosemisimple, 
this proposition is trivial. Therefore, it was not spelled out explicitly 
in our previous papers, where we dealt exclusively 
with the semisimple case.} $\Box$
\end{remark}

This proposition implies that there exists an isomorphism of coalgebras 
$\phi: A\to {\cal O}(G)$ that induces $F.$
Therefore, we can naturally identify the vector spaces 
$V$ and $F(V)$ for all $V\in {\cal C}$, in a functorial way. 

Now, recall that $F$ has a tensor
structure. Namely, we have a collection of right \linebreak ${\cal
O}(G)-$comodule
isomorphisms
$J_{VW}:F(V)\otimes F(W)\raro F(V\otimes W)$ indexed by all
pairs $V,W\in {\cal C}.$ But for any $U\in {\cal C},$ we have
already identified the vector spaces $U$ and $F(U).$ Thus, we can
regard the tensor structure as a collection of isomorphisms 
of vector spaces $V\otimes
W\raro V\otimes W,$ which is functorial with respect to $V$
and $W.$
This collection isomorphisms defines an 
element $J\in ({\cal O}(G)\ot {\cal O}(G))^*$
(see the preceding remarks to Lemma \ref{cv}). It can be
checked (similarly to Theorem 
2.1 in [EG1]) that $J$ is a Hopf 2-cocycle, and 
that $\phi$ induces an isomorphism of cotriangular 
Hopf algebras $A\to {\cal O}(G)^J$.
The implication $(i)\Rightarrow (ii)$ is proved. 

$(ii)\Rightarrow (i)$.

We may assume that $A=
{\cal O}(G)^{J}$ (since $S$ does not depend on the cotriangular
structure). 
 Since the categorical dimension
of any object $V\in {\cal C}$ does
not change under twisting, we have $\dim_c(V)=\dim(V)$ in the
category 
of finite-dimensional right comodules over $A$. 
Therefore, we have 
\begin{eqnarray*}
\lefteqn{\tr (S^2|_{A_V})}\\
& = & \sum_{X,Y\in \Irr({\cal C})}
N_V(X,Y)\dim_c(X)\dim_c(Y)\\
& = & 
\sum _{X,Y\in \Irr({\cal C})}N_V(X,Y)\dim(X)\dim(Y)\\
& = & \dim (A_V).
\end{eqnarray*}
Since any finite-dimensional subcoalgebra $C$ of $A$ has the form
$A_V$
(for $V=C$), this completes the
proof of the theorem. \qed
\subsection{Proof of Proposition \ref{morita}}

Let us first prove the proposition in the 
case when $A,B$ are finite-dimensional.

We need the following standard theorem from noncommutative algebra, which 
can be found for example in [DK, Chapter 3].

Let ${\cal A}$ be a finite-dimensional algebra. 
Let $\Irr({\cal A})$ be the set of isomorphism classes of
irreducible left ${\cal A}$-modules. 
For $M\in \Irr({\cal A})$, let $P(M)$ be the projective 
cover of $M$. 
\begin{theorem}\label{proj} (i) Any finite-dimensional projective 
${\cal A}$-module is a direct sum of $P(M)'s.$ 

\noin
(ii) For any $M\in \Irr({\cal A}),$ the multiplicity of $P(M)$ in the
regular representation ${\cal A}$ is equal to $\text{dim}(M).$
\end{theorem} 

Now we prove the proposition (in the finite-dimensional case).
Let $F_A,F_B$ be the forgetful functors 
from the categories of right $A$-comodules and right $B$-comodules,
respectively, to the category of vector spaces.
All we need to show is that $F_B\circ F$ is isomorphic to $F_A$. 
Indeed, we have $\text{End}(F_A)=A^*$ and
$\text{End}(F_B\circ F)=
B^*$,
so any isomorphism between these two functors will induce 
an isomorphism of coalgebras $A\to B$, which (as one can easily see)
induces $F$. 

The functor $F_A$ is represented by the regular representation $A^*$, 
and $F_B\circ F$ by $F^{-1}(B^*)$. So it suffices to prove
that $F^{-1}(B^*)$ is isomorphic to $A^*$ as 
an $A$-comodule. 

Since $B^*$ is free, it is projective, so $F^{-1}(B^*)$ is also
projective
(as projectivity, unlike freeness, is a categorical property). 
Thus, by Theorem \ref{proj}, we have 
$$
A^*=\bigoplus_{M\in \text{Irr}(A^*)}\text{dim}(M)P(M)\;\text {and}\; 
F^{-1}(B^*)=\bigoplus_{M\in \text{Irr}(A^*)}x(M)P(M),
$$ 
where $x(M)$ are nonnegative integers. 
But since $F$ preserves dimensions, 
we have for any $M\in \text{Irr}(A^*)$: 
$$
\text{dim}(M)=\text{dim}(F(M))=\text{dim}\left(\Hom_{B^*}(B^*,F(M)) 
\right)=
\text{dim}\left(\Hom_{A^*}(F^{-1}(B^*),M)\right)=x(M).
$$
This completes the proof of the proposition
in the finite-dimensional case. 

Now let us consider the infinite-dimensional case.
For simplicity consider the case when $A,B$ are 
countably dimensional (the general case is similar). 
Then $A=\cup_{n\ge 1} A_n,$
where $A_n$ are finite-dimensional coalgebras. 
Let $F_n: \text{Comod}_{f.d}(A_n)\to 
\text{Comod}_{f.d}(B)$ be the restriction of $F,$ and 
let $B_n:=\text{End}(\text{Forget}\circ F_n)^*,$ 
where Forget is the forgetful functor on $B$-comodules. 
It is clear that $B=\cup_{n\ge 1} B_n.$

It is clear from the above finite-dimensional proof 
that for any $n,$ there exists an isomorphism of coalgebras
$\phi_n: A_n\to B_n,$ such that $\phi_{n+1}|_{A_{n}}=
\phi_{n}\circ \text{Ad}(a_{n})^*,$ where 
$a_{n}\in A_{n}^*$ is an invertible element
(this follows from the fact that $\phi_n$
comes from an isomorphism of functors). 
Since the map of the multiplicative groups $Gr(A^*)\to
Gr(A_{n}^*)$
is surjective (because so is the corresponding Lie algebra map), 
we can lift $a_{n}$ to an invertible element of $A^*.$ Abusing
notation, we will denote this element also by $a_{n}.$

Define $\psi_n:=\phi_n\circ
\text{Ad}(a_1^{-1}\cdots a_{n-1}^{-1})^*.$ 
Then $\psi_{n+1}|_{A_{n}}=\psi_{n}$, so $\lbrace \psi_n\rbrace$
defines an isomorphism of coalgebras $\psi: A\to B$ which 
induces a functor isomorphic to $F.$ \qed
\section{Examples of Twisted Function Algebras}
In this section we give examples of Hopf $2-$cocycles
for ${\cal O}(G)$ for certain algebraic groups $G$.
We will construct these cocycles in the form 
of linear endomorphisms of tensor products of any two
finite-dimensional
$G$-modules $V,W$, functorial in terms of $V,W$
(cf. the discussion before Lemma \ref{cv}).  
These examples can be generalized, as usual, using the fact that 
if $J$ is a Hopf 2-cocycle for $G$ and $\phi:G\to G'$ is a
homomorphism then 
$(\phi\otimes \phi)(J)$ is a Hopf 2-cocycle for 
$G'$.
\begin{Example}\label{1} {\rm Let $G$ be the group of
translations of an affine space.
Let $\g$ be the Lie algebra of $G.$ Clearly, $\g$ is abelian.
Therefore, it is straightforward to verify that for any $r\in
\Lambda^2\g,$ the element $J(h):=e ^{hr/2}$ is a Hopf
2-cocycle for any $h\in k$ (the exponential series terminates
in any $V\otimes W$ as the components of $r$ are nilpotent,
so we get a polynomial of $h$). In this case 
$u=1$ and $S^2$ is the identity.} $\Box$
\end{Example}
\begin{Example}\label{2}
{\rm Let $G$ be the group of affine transformations of the line.
Its Lie algebra $\g$ is spanned by two elements 
$X,Y$ such that $[X,Y]=Y$. Define 
$$
J(h):=\sum_{n\ge 0} \frac{h^n}
{n!} X(X-1)\cdots(X-n+1)\ot Y^n.
$$ 
It is not difficult to check that this is a
Hopf 2-cocycle. In this case $u=1+hY+O(h^2),$ so $S^2$ is
not the identity. Thus ${\cal O}(G)^J$ is an example of 
a pseudoinvolutive but not
involutive cortiangular Hopf algebra. 
Such a Hopf algebra can even be cosemisimple: 
it is enough to naturally embed $G$ into $GL(2)$ and 
consider the cosemisimple Hopf algebra ${\cal O}(GL_2)^J$.} $\Box$
\end{Example}
\begin{Example}\label{3}
{\rm For every nonnegative integer $n,$ let $S_n$ be the symmetric
group
of permutations of $n$ symbols. For any $s\in S_n$ and $0\le
m\le n,$
and any Lie algebra $\g,$ let 
$L_{s,m}: \g^{\ot n}\raro U(\g)^{\ot 2}$ be the linear map
determined by 
$$L_{s,m}(a_1\ot\cdots\ot a_n)=
a_{s(1)}\cdots a_{s(m)}\ot a_{s(m+1)}\cdots a_{s(n)}.$$
Let $X$ be the disjoint union of the sets $X_n:=S_{2n}\times
\{0,...,2n\}$
for $n\ge 2.$ Let $k[X]$ be the set of $k-$valued functions
on $X.$ For any $f\in k[X],$ and any Lie algebra $\g,$
define the function
$J_f:\Lambda^2 \g\raro U(\g)^{\ot 2}[[h]],$ by the formula
\begin{equation}\label{jfrh}
J_f(r,h)=1+hr/2+\sum_{n\ge 2}h^n\sum_{(s,m)\in
X_{n}}f(s,m)L_{s,m}(r^{\ot n}).
\end{equation}
}
\end{Example}

It is easy to show that if $J_f(r,h)$ is a twist for $U(\g)[[h]],$
then $r$ must satisfy the classical Yang-Baxter equation (CYBE):
\begin{equation}\label{cybe}
[r_{12},r_{13}]+[r_{12},r_{23}]+[r_{13},r_{23}]=0.
\end{equation}
Conversely, let $r$ be a solution of CYBE.
\begin{Definition}\label{uq} 
We say that $f\in k[X]$ is a quantization function for $r$ 
if the element $J_f(r,h)$ is a twist.
We say that $f$ is a universal quantization function 
if this is the case for all $\g,r$.   
\end{Definition}

It is easy to construct concrete examples of quantization
functions. For instance, in Example \ref{1}, it is
straightforward to verify that $e^{hr/2}$ comes
from a quantization function for any $r$ and abelian $\g.$ The
existence of universal quantization functions is not obvious, 
but it follows from the quantization theory of
[EK1]. Namely, 
a construction of such functions can be obtained
from formula (3.1) in [EK1].

Quantization functions allow one to construct 
a large family of examples of Hopf 2-cocycles. Namely, we have:
\begin{Theorem}\label{csag} 
Suppose that $G$ is an algebraic group with Lie algebra $\g,$
and let $N$ be its unipotent radical with Lie algebra $\n.$
Suppose
that $r$
is an
element of $\g\wedge \n \subset \Lambda^2 \g$ which
satisfies the CYBE. Then
for any $f\in k[X],$ and any two algebraic representations $V,W$ of $G,$
$J_f(r,h)|_{V\ot W}$ is a polynomial in $h$ (i.e. the series terminates).
In particular, for any $h\in k,$ $J_f(r,h)$ is a well defined
element in $({\cal O}(G)\ot {\cal O}(G))^*.$ This element is a Hopf
2-cocycle
for ${\cal O}(G)$ if $f$ is a quantization function for $r$. 
\end{Theorem}
\proof We only need to show that the series $J_f(r,h)$ terminates. 
Let $V,W$ be two finite-dimensional algebraic representations of
$G$. Let 
$B_V,B_W$ be the images of $U(\g)$ in
$\text{End}(V)$ and $\text{End}(W),$ and let
$I_V,I_W$ be the nilpotent radicals of $B_V,B_W$. 
It is clear that under the action of $G$ in $V,W$, 
$\n$ maps into $I_V,I_W$. Let $n$ be a positive 
integer such that $I_V^n=I_W^n=0$. 
Then it is clear from the definition 
that the series $J_f(r,h)|_{V\otimes W}$ is a polynomial in $h$ 
of degree $\le 2n-2$. \qed
\section{The Unipotency Conjecture}
We conclude the paper with the following
\begin{Conjecture}\label{unip} 
For any pro-algebraic group $G$ and Hopf $2-$cocycle $J$ for ${\cal
O}(G)$ over $k,$ the operator $S^2$ is unipotent on ${\cal O}(G)^J.$
\end{Conjecture}
\begin{Remark}
{\rm We know that the sum of the eigenvalues of $S^2$ on any 
finite-dimensional subcoalgebra is equal to its dimension, 
but the conjecture says that furthermore all of these eigenvalues
are $1.$} $\Box$
\end{Remark} 
\begin{Remark}\label{ex} 
{\rm Conjecture \ref{unip} is obviously satisfied in Examples \ref{1}
and \ref{2}. Moreover, it follows from the theorem below 
that it is satisfied in Example \ref{3} as well.} $\Box$
\end{Remark} 

Let $\Sigma$ be an irreducible affine algebraic curve with 
a marked smooth 
point $0$, and let ${\cal O}(\Sigma)$ be the ring of regular
functions on $\Sigma$. The standard example is $\Sigma:=k,$ 
${\cal O}(\Sigma)=k[x].$

Let $J:{\cal O}(G)^{\otimes 2}\to {\cal O}(\Sigma)$ 
be a family of Hopf 2-cocycles for ${\cal O}(G)$ parametrized by $a\in
\Sigma,$ with $J(0)=1.$ In this case, we will say that $J(a),$ for
any $a\in \Sigma,$ is obtained by deformation of $J(0).$ 
\begin{Theorem}\label{unipot} Let $J$ be obtained by deformation of $1$. 
Then Conjecture \ref{unip} holds for the cotriangular Hopf
algebra ${\cal
O}(G)^{J}.$
\end{Theorem}

The rest of the section is devoted to the proof of the theorem.

To prove the theorem, we will choose a local parameter 
$h$ on $\Sigma$, and write $J$ as a formal power series
in $h$: $J=1+\sum_{n\ge 1}h^nr_n$,
$r_n\in ({\cal O}(G)^{\otimes 2})^*$. 
We will say that $J$ is {\em local} if $r_n\in U(\g)^{\otimes 2}$ for
all $n$. 
\begin{Lemma}\label{l64} 
The Hopf $2-$cocycle $J$ is gauge equivalent to a local Hopf
$2-$cocycle. That is,
there exists a ``gauge transformation'' $g:=1+hg_1+h^2g_2+...$, 
$g_i\in {\cal O}(G)^*$, $\varepsilon(g_i)=0$, such that 
the Hopf $2-$cocycle $J^g:=\Delta(g)J(g^{-1}\otimes g^{-1})$ is local. 
\end{Lemma} 
\proof 
Let us prove the statement modulo $h^{n+1}$ by induction in $n.$
The base of induction ($n=0$) is clear. 
To do the inductive step, assume that $J$ 
is local modulo $h^n$. Observe that since $r_n$ 
satisfies the Hopf $2-$cocycle condition it follows that 
$$
r_n^{12}+(\Delta\otimes I)(r_n)-
r_n^{23}-(I\otimes \Delta)(r_n)=f(r_1,...,r_{n-1}),
$$
where $f$ is a polynomial. Thus, 
we have $dr_n\in U(\g)^{\otimes 3}$, where $d$ is the
differential in the Hochschild complex of the coalgebra ${\cal O}(G)^*$
with trivial coefficients. 

It is well known that the embedding of coalgebras
$U(\g)\to {\cal O}(G)^*$ defines an isomorphism of Hochschild 
cohomology of these coalgebras 
with trivial coefficients, and that both cohomology spaces 
are equal to $\Lambda \g$, with the usual grading 
(the fact that the cohomology of $U(\g)$ is $\Lambda \g$ is
discussed for example
in [Dr1, p.1435]). Therefore, 
there exists $r_n'\in U(\g)^{\otimes 2}$ such that 
$dr_n'=dr_n$. Let $s:=r_n-r_n'$. Then $ds=0$. 
Therefore, we have $s=s_0+dz=s_0+\Delta(z)-z\otimes 1-1\otimes
z$, 
where $s_0\in \Lambda^2\g$ and $z\in {\cal O}(G)^*$. 
Let us replace $J$ with $J^g$ for $g:=1+h^n(z-\epsilon(z))$. 
Then $J^g$ is local modulo $h^{n+1}$. 
The lemma is proved. 
\qed

Now let us continue the proof of the theorem.
By Lemma \ref{l64}, we can assume that $J(h)$ is local. 
Then $J(h)$ is a twist for $U(\g)[[h]]$, so one can 
define the triangular quantized universal enveloping (QUE) 
algebra $U(\g)[[h]]^{J(h)}$. 
It is sufficient to show that the Drinfeld element $u$ of 
this QUE algebra is unipotent on every finite-dimensional 
representation of $\g$.

Using the main result of [EK2] (in the triangular case), 
we conclude that $U(\g)[[h]]^{J(h)}$ is isomorphic to $U_h(\g,r)$, where
$U_h$ is the quantization functor form [EK2], and 
$r\in \Lambda^2\g[[h]]$ is a solution of CYBE.

The QUE algebra $U_h(\g,r)$, by definition, is obtained by 
twisting $U(\g)$ using a twist $J_f(r,h)$, where $f$ 
is a universal quantization function. This implies that 
the element $u$
in this QUE algebra has the form 
$$
u=1+\sum_{n\ge 1}\sum_{\sigma\in S_{2n}}
c_{n,\sigma}m_{2n}(\sigma\circ r^{\otimes n}),
$$ where $m_{2n}$ is the multiplication of $2n$
components.  

Now we show the unipotency of $u$ in any 
finite-dimensional $\g$-module by induction in the rank of
$r$. Without loss of generality, we can assume that 
$\g$ is spanned by the components of $r$, and that $V$ is 
an irreducible faithful $\g$-module (if it is not faithful, we
can go to r-matrices of smaller rank, for which the statement is
known). This implies that $\g$ is reductive. 

Now we will use the following theorem of Drinfeld (see e.g. 
[ES, Proposition 5.2]): 
\begin{Theorem}\label{drin} Let $\g$ be a Lie algebra and 
$r\in \Lambda^2\g$ be a solution of CYBE whose components span 
$\g$. If $\g$ is reductive, then it is abelian.
\end{Theorem}

This theorem immediately implies Theorem \ref{unipot}: if $\g$ is
abelian then, because of the skew-symmetry of $r$, 
we have $u=1$. 

The rest of this section is the proof of Theorem \ref{drin}, 
which we give for the reader's convenience. 

It is clear that $r$ defines a nondegenerate 
skew-symmetric bilinear form on $\g^*$, hence one can define 
a skew-symmetric form $r^{-1}$ on $\g$. 
Since $r$ satisfies CYBE, this form is a 2-cocycle. 

Let $\g=\g_s\oplus \g_a$ be the splitting of $\g$ into the
semisimple and the abelian parts. Assume that $\g_s\ne 0$.  
Then $H^2(\g,k)=H^2(\g_a,k)=\Lambda^2\g_a^*$, which shows that 
$r^{-1}$ can be written as 
$r^{-1}(x,y)=f([x,y])+\rho(x,y)$, where $\rho\in \Lambda^2\g_a^*$,
$f\in \g^*$.  
Thus, the decomposition $\g=\g_s\oplus \g_a$ is orthogonal 
under $r^{-1}$, and hence the form $f([x,y])$ 
has to be nondegenerate on $\g_s$. 

But $\g_s$ has a nondegenerate invariant form $(,)$, and we can
write the functional 
$f(x)$ as $(z,x)$ for some $z\in \g_s$. 
Thus, the form $(z,[x,y])$ should be nondegenerate. 
However, this is impossible as $z$ obviously belongs 
to the kernel of this form (since $(z,[z,y])=0$).  Thus, we have a
contradiction and 
hence $\g_s=0$. 
\qed

\section{Questions}
In conclusion let us discuss possible directions for future
research. 
The main remaining problem is to obtain a classification of
twisted group algebras. We believe that this should be done 
by generalizing the techniques of 
[M] and [EG2] to the infinite-dimensional case and 
combining them with the techniques of [EK1,EK2].   
Let us formulate some precise questions which are 
related to this problem. 
\begin{Question} {\rm Is any Hopf 2-cocycle for ${\cal
O}(G)$ gauge equivalent to a deformation of a Hopf 2-cocycle
of finite rank?}
\end{Question}
\begin{Question}\label{q2}{\rm 
Say that a cotriangular Hopf algebra is {\em minimal} if the R-form
is nondegenerate. 
Now suppose that $A$ is a minimal pseudoinvolutive cotriangular
Hopf algebra, and that the underlying group $G$ is reductive.  
Is it true that the connected component of the identity in
$G$ is abelian (i.e. a torus)?}
\end{Question}
\begin{Remark} {\rm The answer to the classical analog of
Question \ref{q2} is positive: if a Lie algebra spanned
by the components of a skew-symmetric solution of CYBE is
reductive, then it is abelian (see Theorem \ref{drin}).} $\Box$ 
\end{Remark}
\begin{Question}\label{q3}{\rm
Let $A$ be any cotriangular Hopf algebra over $k.$ Is it true that
the eigenvalues of the square of its antipode are all roots of unity?
Are all $\pm 1$? Is it true at least for twisted function algebras?}
\end{Question}

Note that a positive answer to either form of 
Question \ref{q3} will imply Conjecture
\ref{unip}.

\end{document}